\documentclass[12pt]{article}
\usepackage{amscd}
\usepackage{amsmath,amsfonts,amssymb,amscd}
\usepackage{indentfirst,graphics,epsfig,psfrag}
\usepackage{ifpdf}
\input{epsf}

\makeatletter \@addtoreset{figure}{section} \makeatother
\makeatletter
\long\def\@makecaption#1#2{%
   \vskip 10\p@
   \setbox\@tempboxa\hbox{{#1}\ \ #2}%
   \ifdim \wd\@tempboxa >\hsize
       {#1}\ \ #2\par
   \else
       \hbox to\hsize{\hfil\box\@tempboxa\hfil}%
   \fi}
\makeatother

\input{epsf}

\newtheorem{thm}{Theorem}

\newtheorem{lem}{Lemma}

\begin{document}
\rule{0cm}{1cm}
\begin{center}
{\Large\bf Nordhaus-Gaddum-type theorem for the\\[2mm] rainbow
vertex-connection number of a graph\footnote{Supported by NSFC No.
11071130.}}
\end{center}

\begin{center}
Lily Chen, Xueliang Li, Mengmeng Liu\\
Center for Combinatorics, LPMC\\
Nankai University, Tianjin 300071, P. R. China\\
Email:  lily60612@126.com, lxl@nankai.edu.cn, liumm05@163.com
\end{center}

\begin{center}
\begin{minipage}{120mm}
\vskip 0.3cm

\begin{center}
{\bf Abstract}
\end{center}

{\small  A vertex-colored graph $G$ is rainbow vertex-connected if
any pair of distinct vertices are connected by a path whose internal
vertices have distinct colors. The rainbow vertex-connection number
of $G$, denoted by $rvc(G)$, is the minimum number of colors that
are needed to make $G$ rainbow vertex-connected. In this paper we
give a Nordhaus-Gaddum-type result of the rainbow vertex-connection
number. We prove that when $G$ and $\overline{G}$ are both
connected, then $2\leq rvc(G)+rvc(\overline{G})\leq n-1$. Examples
are given to show that both the upper bound and the lower bound
are best possible for all $n\geq 5$. }\\[3mm]
\end{minipage}
\end{center}

\section{Introduction}

All graphs considered in this paper are simple, finite and
undirected. We follow the notation and terminology of
\cite{Bondy-Murty}. An edge-colored graph $G$ is rainbow connected
if any pair of distinct vertices are connected by a path whose edges
have distinct colors. Clearly, if a graph is rainbow edge-connected,
then it is also connected. Conversely, any connected graph has
trivial edge coloring that makes it rainbow edge-connected; just
color each edge with a distinct color. The rainbow connection number
of a connected graph $G$, denoted by $rc(G)$, is the minimum number
of colors that are needed in order to make $G$ rainbow connected,
which was introduced by Charrand et al. Obviously, we always have
$diam(G)\leq rc(G)\leq n-1$, where $diam(G)$ denotes the diameter of
$G$. Notice that $rc(G)=1$ if and only if $G$ is a complete graph,
and that $rc(G)=n-1$ if and only if $G$ is a tree.

In \cite{KY}, Krivelevich and Yuster proposed the concept of rainbow
vertex-connection. A vertex-colored graph is rainbow
vertex-connected if any pair of distinct vertices are connected by a
path whose internal vertices have distinct colors. The rainbow
vertex-connection of a connected graph $G$, denoted by $rvc(G)$, is
the minimum number of colors that are needed to make $G$ rainbow
vertex-connected. An easy observation is that if $G$ is a connected
graph with $n$ vertices then $rvc(G)\leq n-2$. We note the trivial
fact that $rvc(G)=0$ if and only if $G$ is a complete graph. Also,
clearly, $rvc(G)\geq diam(G)-1$ with equality if the diameter is $1$
or $2$.

A Nordhaus--Gaddum-type result is a (tight) lower or upper bound on
the sum or product of a parameter of a graph and its complement. The
name Nordhaus--Gaddum-type is used because in 1956 Nordhaus and
Gaddum \cite{Nordhaus-Gaddum} first established the following
inequalities for the chromatic numbers of graphs, they proved that
if $G$ and $\overline{G}$ are complementary graphs on $n$ vertices
whose chromatic numbers are $\chi(G)$, $\chi(\overline{G})$
respectively, then
$$2\sqrt{n}\leq \chi(G)+\chi(\overline{G}) \leq n+1.$$ Since then, many
analogous inequalities of other graph parameters are concerned, such
as domination number \cite{Harary-Haynes}, Wiener index and some
other chemical indices \cite{Zhang-Wu}, and so on.

In \cite{LCL}, the authors considered Nordhaus--Gaddum-type result
for the rainbow connection number. In this paper, we are concerned
with analogous inequalities involving the rainbow vertex-connection
number of graphs. We prove that $$2\leq rvc(G)+rvc(\overline{G})\leq
n-1.$$

The rest of this paper is organized as follows. Section 2 contains
the proof of the sharp upper bound. Section 3 contains the proof of
the sharp lower bound.

\section{Upper bound on $rvc(G)+rvc(\overline{G})$}

We begin this section with two lemmas that are needed in order to
establish the proof of the upper bound.

\begin{lem}\label{thm1}
Let $G$ be a nontrivial connected graph of order $n$, and
$rvc(G)=k$. Add a new vertex $v$ to $G$, and make $v$ be adjacent to
$q$ vertices of $G$, the resulting graph is denoted by $G'$. Then if
$q\geq n-k$, we have $rvc(G')\leq k.$
\end{lem}

\begin{pf}Let $c: V(G) \rightarrow \{1, 2, \cdots, k\}$ be a rainbow
$k$-vertex-coloring  of $G$, $X=\{x_1, x_2, \cdots, x_q\}$ be the
vertices that are adjacent to $v$, $V\backslash X=\{y_1, y_2,
\cdots, y_{n-q}\}$. We can assume that there exists some $y_j$ such
that there is no rainbow vertex-connected-path from $v$ to $y_{j}$;
otherwise, the result holds obviously. Because $G$ is a rainbow
$k$-vertex-coloring, there is a rainbow vertex-connected-path $P_i$
from $x_i$ to $y_j$ for every $x_i$, $i\in \{1,2,\cdots, q\}$.
Certainly, $P_i \bigcap P_j$ may not be empty. We notice that no
other vertices of $\{x_1,x_2,\cdots,x_q\}$ different from $x_i$
belong to $P_i$ for each $1\leq i \leq q$. If so, let ${x_i}'$ be
the last vertex in $\{x_1,x_2,\cdots,x_q\}$ which belongs to $P_i$,
denote $P_i$ by $x_i{P_i}'{x_i}'Q_iy_j$, then $v{x_i}'Q_iy_j$ is a
rainbow vertex-connected-path, a contradiction to our assumption.
Since $v$ and $y_j$ are not rainbow vertex-connected, for each
$P_i$, there is some $y_{k_i}$ such that $c(x_i)=c(y_{k_i})$. That
means that the colors that are assigned to $X$ are among the colors
that are assigned to $V\backslash X $. So $rvc(G)=k\leq n-q$. By the
hypothesis $q\geq n-k$, we have $rvc(G)= n-q$, that is, all vertices
in $V\backslash X$ have distinct colors. Now we construct a new
graph $G'=P_1\bigcup P_2\bigcup \cdots \bigcup P_q$. For every $y_t$
not in $G'$, there is a $y_s \in G'$ such that $y_ty_s \in E(G)$. Or
$N(y_t)\subseteqq \{x_1,x_2,\cdots, x_q\}$. Since $G$ is rainbow
$k$-vertex-connected, there is a rainbow vertex-connected path from
$y_t$ to $y_j$, denoted by $y_tx_kQy_j$, where $x_k \in N(y_t)$.
Thus $vx_kQy_j$ is a rainbow vertex-connected path, a contradiction.
It follows that $G[y_1,y_2,\cdots,y_{n-q}]$ is connected. Certainly,
$G[y_1,y_2,\cdots,y_{n-q}]$ has a spanning tree $T$, and $T$ has at
least two pendant vertices. Then there must exist a pendant vertex
whose color is different from $x_1$, and we assign the color to
$x_1$. It is easy to check that $G$ is still rainbow
$k$-vertex-connected, and there is a rainbow vertex-connected path
between $v$ and $y_j$. If there still exists some $y_j$ such that
$v$ and $y_j$ are not rainbow vertex-connected, we do the same
operation, until $v$ and $y_j$ are rainbow vertex-connected for each
$j\in\{1,2,\cdots, n-q\}$. Thus $G'$ is rainbow vertex-connected. It
follows that $rvc(G')\leq k.$

\end{pf}
\begin{qed}
\end{qed}

\begin{lem}\label{lem1}
Let $G$ be a connected graph of order $5$. If $\overline{G}$ is
connected, then $rvc(G)+rvc(\overline{G})\leq 4$.
\end{lem}

\begin{pf}We consider the situations of $G$.

First, if $G$ is a path, then $rvc(G)=3$. In this case
$diam(\overline{G})=2$, and then $rvc(\overline{G})=1$.

Second, if $G$ is a tree but not a path, then $rvc(G)<3$. Since $G$
is a bipartite graph, then $\overline{G}$ consists of a $K_2$ and a
$K_3$ and two edges between them. So we assign color $1$ to the
vertices of $K_2$ and color $2$ to the vertices of $K_3$, and this
makes $\overline{G}$ rainbow vertex-connected, that is,
$rvc(\overline{G})\leq 2$.

Finally, if both $G$ and $\overline{G}$ are not trees, then
$e(G)=e(\overline{G})=5$. If $G$ contains a cycle of length $5$,
then $G=\overline G=C_5$, thus $rvc(G)=rvc(\overline{G})=1$. If $G$
contains a cycle of length $4$, there is only one graph $G$ which is
showed in Figure 1, we can color $G$ and $\overline G$ with 2 colors
to make them rainbow vertex-connected, see Figure 1. If $G$ contains
a cycle of length $3$, then $G$ and $\overline G$ are showed in
Figure 2. By the coloring showed in the graphs, we have
$rvc(G)+rvc(\overline{G})=4.$

\setlength{\unitlength}{1.2mm}
\begin{center}
\begin{picture}(80,60)
\put(0,45){\circle*{1}} \put(10,45){\circle*{1}}
\put(0,55){\circle*{1}} \put(10,55){\circle*{1}}
\put(20,55){\circle*{1}} \put(50,45){\circle*{1}}
\put(50,55){\circle*{1}} \put(60,50){\circle*{1}}
\put(70,50){\circle*{1}} \put(80,50){\circle*{1}}

\put(-0.1,45){\line(0,1){10}} \put(0,45){\line(1,0){10}}
\put(0,55){\line(1,0){20}}\put(10,45){\line(0,1){10}}
\put(50,45){\line(0,1){10}}\put(60,50){\line(1,0){20}}
\put(60,50){\line(-2,1){10}}\put(50,45){\line(2,1){10}}

\put(5,13){\circle*{1}} \put(15,13){\circle*{1}}
\put(5,20){\circle*{1}} \put(15,20){\circle*{1}}
\put(10,27){\circle*{1}} \put(55,13){\circle*{1}}
\put(65,13){\circle*{1}} \put(55,20){\circle*{1}}
\put(65,20){\circle*{1}} \put(60,27){\circle*{1}}

\put(5,13){\line(0,1){7}} \put(15,13){\line(0,1){7}}
\put(5,20){\line(1,0){10}}\put(5,20){\line(3,4){5}}
\put(15,20){\line(-3,4){5}}\put(55,13){\line(0,1){7}}
\put(65,13){\line(0,1){7}}
\put(55,20){\line(1,0){10}}\put(55,20){\line(3,4){5}}
\put(65,20){\line(-3,4){5}}

\put(-0.5,56){\tiny{$1$}}\put(10,56){\tiny{$2$}}
\put(60,51){\tiny{$1$}}\put(70,51){\tiny{$2$}}

\put(4.5,21){\tiny{$1$}}\put(15,21){\tiny{$2$}}
\put(54.5,21){\tiny{$1$}}\put(65,21){\tiny{$2$}}

\put(10,40){\small$G$}\put(60,40){\small$\overline{G}$}
\put(10,32){Figure $1$: $G$ contains a cycle of length 4.}

\put(10,8){\small$G$}\put(60,8){\small$\overline{G}$}
\put(8,2){Figure $2$: $G$ contains a cycle of length 3.}

\end{picture}
\end{center}

By these cases, we have $rvc(G)+rvc(\overline{G})\leq 4.$
\end{pf}
\begin{qed}
\end{qed}

From the above lemmas, we have our first theorem.
\begin{thm}
$rvc(G)+rvc(\overline{G})\leq n-1$ for all $n\geq 5$, and this bound
is best possible.
\end{thm}

\begin{pf}
We use induction on $n$. By Lemma \ref{lem1}, the result is evident
for $n=5$. We assume that $rvc(G)+rvc(\overline{G})\leq n-1$ holds
for complementary graphs on $n$ vertices. To the union of a
connected graph $G$ and its $\overline{G}$, which forms the complete
graph on these $n$ vertices, we adjoin a vertex $v$. Let $q$ of the
$n$ edges between $v$ and the union be adjoined to $G$ and the
remaining $n-q$ edges to $\overline{G}$. If $G'$ and $\overline{G'}$
are the graphs so determined (each of order $n+1$), then
$$rvc(G')\leq rvc(G)+1, \ \ rvc(\overline{G'})\leq rvc(\overline{G})+1.$$
These inequalities are evident from the fact that if given a rainbow
$rvc(G)$-vertex-coloring ($rvc(\overline{G})$-vertex-coloring) of
$G$ ($\overline{G}$), we assign the new color to the vertex which is
adjacent to $v$ and keep other vertices unchanged, the resulting
coloring makes $G'$ ($\overline{G'}$) rainbow vertex-connected. Then
$rvc(G')+rvc(\overline{G'})\leq rvc(G)+rvc(\overline{G})+2\leq n+1.$
And $rvc(G')+rvc(\overline{G'})\leq n$ except possibly when
$$rvc(G')=rvc(G)+1, \ \ rvc(\overline{G'})=rvc(\overline{G})+1.$$
In this case, by Lemma \ref{thm1}, $q\leq n-rvc(G)-1, n-q\leq
n-rvc(\overline{G})-1$, thus $rvc(G)+rvc(\overline{G})\leq n-2,$
from which $rvc(G')+rvc(\overline{G'})\leq n$. This completes the
induction.

The following example shows that the bound established is sharp for
all $n \geq 5$: If $G$ be a path of order $n$, then $rvc(G)= n-2$.
It is easy to obtain $\overline{G}$, and check that
$diam(\overline{G})=2$. Then $rvc(\overline{G})= 1$, and so we have
$rvc(G)+rvc(\overline{G})= n-1$.
\end{pf}
\begin{qed}
\end{qed}

\section{Lower bound on $rvc(G)+rvc(\overline{G})$}

As we note that $rvc(G)=0$ if and only if $G$ is a complete graph.
Thus if we want both $G$ and $\overline{G}$ are connected, and so
$rvc(G)\geq 1, rvc(\overline{G})\geq 1$. Then $rvc(G)+
rvc(\overline{G})\geq 2.$ Our next theorem shows that the lower
bound is sharp for all $n \geq 5$.
\begin{thm}\label{thm2}
For $n \geq 5$, the lower bound of $rvc(G)+ rvc(\overline{G})\geq 2$
is best possible, that is, there are graphs $G$ and $\overline{G}$
with $n$ vertices, such that $rvc(G)= rvc(\overline{G})=1.$
\end{thm}

\begin{pf}We only need to prove that for $n\geq 5$, there are graphs $G$ and
$\overline{G}$ with $n$ vertices, such that $diam(G)=
diam(\overline{G})=2.$

We construct $G$ as follows: if $n=2k+1$,
$$
V(G)=\{v,v_1,v_2,\cdots,v_k,u_1,u_2,\cdots,u_k\}
$$
$$
E(G)=\{vv_i|1\leq i\leq k\}\bigcup \{v_iu_i|1\leq i\leq k\}\bigcup
\{u_iu_j|1\leq i,j\leq k\};
$$
if $n=2k$,
$$
V(G)=\{v,v_1,v_2,\cdots,v_k,u_1,u_2,\cdots,u_{k-1}\}
$$
$$
E(G)=\{vv_i|1\leq i\leq k\}\bigcup \{v_iu_i|1\leq
i<k\}\bigcup\{v_ku_{k-1}\}\bigcup \{u_iu_j|1\leq i,j\leq k-1\}.
$$
We can easily check that $diam(G)= diam(\overline{G})=2$.
\end{pf}
\begin{qed}
\end{qed}


\begin{thebibliography}{99}

\bibitem{Bondy-Murty} J.A. Bondy, U.S.R. Murty,
Graph Theory, GTM 244, Springer, 2008.

\bibitem{Chartrand-Johns} G. Chartrand, G.L. Johns, K.A. McKeon, and P. Zhang,
Rainbow connection in graphs, Math. Bohem. 133(2008), 85-98.

\bibitem{KY}
M. Krivelevich, R. Yuster, The rainbow connection of a graph is (at
most) reciprocal to its minimum degree, J. Graph Theory 63(2010),
185--19

\bibitem{Nordhaus-Gaddum} E.A. Nordhaus, J.W. Gaddum, On complementary graphs, Amer.
Math. Monthly 63(1956), 175-177.

\bibitem{Harary-Robinson} F. Harary and R.W. Robinson, The diameter of a graph and its
complement, The American Mathematical Monthly 92(1985), 211-212.

\bibitem{Harary-Haynes} F. Harary and T.W. Haynes, Nordhaus-Gaddum inequalities
for domination in graphs, Discrete Math. 155(1996), 99-10.

\bibitem{Zhang-Wu} L. Zhang and B. Wu, The Nordhaus-Gaddum type inequalities of
some chemical indices, MATCH Commun. Math. Comput. Chem. 54(2005),
189-194.

\bibitem{LCL}
L. Chen, X. Li and H. Lian, Nordhaus-Gaddum-type theorem for rainbow
connection number of graphs, arXiv:1012.2641.

\end{thebibliography}
\end{document}